\documentclass[conference]{IEEEtran}
\usepackage{cleveref}
\usepackage{ifpdf}

\usepackage[noadjust]{cite}

\ifCLASSINFOpdf
  \usepackage[pdftex]{graphicx}
   \graphicspath{{../images/},{../../SIMU/plots/keep/},{../../SIMU/plots/keep/smartGridComm/},{.}}
   \DeclareGraphicsExtensions{.pdf,.jpeg,.png,.svg}
\else
   \usepackage[dvips]{graphicx}
   \graphicspath{{../eps/}}
   \DeclareGraphicsExtensions{.eps}
\fi
\usepackage{amsmath}
\usepackage{algorithmic}
\usepackage{bm}
\usepackage{textcomp}
\usepackage{subfig}

\newcommand{\norm}[1]{\left\Vert #1\right\Vert}

\newcommand{\SC}{\text{S\hspace{-1pt}C}}

\newcommand{\eqDef}{:=}

\newtheorem{theorem}{Theorem}
\newtheorem{Cor}{Corollary}

\newtheorem{definition}{Definition}

\Crefname{equation}{Eq.}{Eqs.}
\Crefname{figure}{Fig.}{Figs.}
\Crefname{tabular}{Tab.}{Tabs.}
\Crefname{table}{Tab.}{Tabs.}
\Crefname{theorem}{Thm.}{Thms.}
\Crefname{definition}{Def.}{Defs.}
\Crefname{section}{Sec.}{Secs.}

\usepackage[colorinlistoftodos,bordercolor=orange,backgroundcolor=orange!20,linecolor=orange,textsize=scriptsize]{todonotes}
\usepackage{multirow}
\usepackage{pbox}
\usepackage{amsfonts}
\newcommand{\mr}{\mathrm}

\newcommand{\x}{\ell}

\newcommand{\tsfrac}[2]{\textstyle\frac{#1}{#2}}
\newcommand{\tssum}{\textstyle\sum}
\begin{document}

\title{Demand Response in the Smart Grid: the Impact of Consumers Temporal Preferences}

\author{\IEEEauthorblockN{Paulin Jacquot\IEEEauthorrefmark{1}\IEEEauthorrefmark{2} \IEEEmembership{Student Member, IEEE},
Olivier Beaude\IEEEauthorrefmark{1}, 
St\'ephane Gaubert\IEEEauthorrefmark{2}, Nadia Oudjane\IEEEauthorrefmark{1}}
\IEEEauthorblockA{\IEEEauthorrefmark{1} EDF, Saclay, France\\
}
\IEEEauthorblockA{\IEEEauthorrefmark{2} Inria Saclay, CMAP, Ecole Polytechnique, Palaiseau, France  }
\\ \textit{Accepted for presentation at IEEE Int. Conf. on Smart Grid Communications (SmartGridComm), 2017, Dresden, Germany.}

}

\maketitle


\begin{abstract}
In Demand Response programs, price incentives might not be sufficient to modify residential consumers load profile. Here, we consider that each consumer has a preferred profile and a discomfort cost when deviating from it. Consumers can value this discomfort at a varying level that we take as a parameter. This work analyses Demand Response as a game theoretic environment. We study the equilibria of the game between consumers with preferences within two different dynamic pricing mechanisms, respectively the daily proportional mechanism introduced by Mohsenian-Rad \textit{et al}, and an hourly proportional mechanism. We give new results about equilibria as functions of the preference level in the case of quadratic system costs and prove that, whatever the preference level, system costs are smaller with the hourly mechanism. We simulate the Demand Response environment using real consumption data from \textit{PecanStreet} database. While the Price of Anarchy remains always close to one up to 10$^{-3}$ with the hourly mechanism, it can be more than 10\% bigger with the daily mechanism.
\end{abstract}


%
\IEEEpeerreviewmaketitle

\section{Introduction}

On the path to a ``grid of the future" \cite{ipakchi2009}, Demand Side Management (DSM) and Demand Response (DR) programs are envisioned to provide a substantial support to face the challenges to come: integration of new  and additional electrical usages, e.g. electric vehicles \cite{beaude2016reducing}, transition to a significant penetration rate of Renewable Energy Sources \cite{wu2012} and a decentralized energy production and reduction of carbon emissions.
Various models and mechanisms have been proposed in the Smart Grid literature to operate the flexibility offered to the system by DSM \cite{deng2015}. Among them, we focus in this paper in the line of game theoretic schemes, introduced with the seminal paper \cite{mohsenian2010autonomous}. In this case, an ``energy consumption game" is defined between consumers --- considered as players. With this interpretation, an Independent System Operator has the responsability to design good rules to induce an efficient collective behaviour of flexible consumers (typically, at Nash equilibria), aligned with electricity system objectives. This efficiency is often measured with system metrics, e.g.  generating or distributing electricity at a reasonable cost / emission level, limited impact on electrical assets, etc. This efficiency can also integrate consumer metrics, which leads to considering social cost as the sum of system cost and consumers objectives.

Individual and personal constraints will influence the behaviour of flexible consumers. 
In turn, many papers have adressed the question of modeling  individual effort of consumers in the context of DSM. Several works relate this individual effort to the desired indoor temperature. For instance, \cite{nguyen2014} considers a distance between ``desired indoor temperature" profile and effective one, weighted by an occupancy variable
. Another standard model consists in penalizing the delay between possible operation time (e.g., the starting period of availability) and effective one of a flexible electrical appliance \cite{yaagoubi2013,chen2011,mohsenian2010opt}. In this case, the cost is generally linear with the waiting time (and sometimes weighted by the power of considered task as in \cite{mohsenian2010opt}). A different approach in \cite{baniya2014} considers as a metric for uncomfort of residential consumers the colour quality of a ``smart lighting".  Note that in all of the aforementioned metrics, the total flexible energy consumed is fixed, and consumption flexibility consists only in temporal scheduling of this fixed amount. Individual effort made by consumers consists then of a temporal preference for consumption. As a result, these metrics can be formulated as particular cases of the framework that we will propose in this paper. 

Other works such as  \cite{li2011optimal,fahrioglu2001using,samadi2012advanced}  rather consider an individual utility term that depends on the total flexible energy consumed; a standard representation of this utility is made with an increasing and concave function of total amount of flexible energy (a quadratic function with a saturation threshold is often used, as mentioned in \cite{deng2015}). In this model, consumers can receive no energy at all and it is assumed that their satisfaction increases with the volume of energy they consume.

Whatever the metrics considered for the individual consumer preference, to the best of our knowledge no study has been made on the impact of the weight given to this preference on the consumers behaviour. More precisely, this weight will influence the induced equilibrium in the associated consumption energy game, which will impact the system efficiency. This is precisely the issue addressed in this paper.

In this work, we will distinguish standard metrics of efficiency in a game: the system cost and the social cost. We will study as well the Price of Anarchy (PoA) \cite{koutsoupias1999worst}, a standard measure of efficiency in a game, and a measure called \emph{Price of Efficiency} defined to be similar to the Price of Anarchy on the system operator side. While the study of such indicators in energy consumption games has been done previously (see \cite{mohsenian2010autonomous} or \cite{beaude2016reducing} which exhibit games where $\mr{PoA}=1$)
, the analysis of the evolution of these indicators with respect to the weight on individual effort term is a novelty adressed here.

This paper brings several contributions. We extend the standard model of an energy consumption game among consumers, studied in \cite{Paulin2017ISGT}, by adding individual temporal preferences. Next, we give theoretical results in this extended framework about the impact of preferences on the equilibria of the game and the efficiency of those equilibria. For that, we analyse the induced social cost and system costs. Last, we present numerical results on a realistic test case, using consumption data from \textit{PecanStreet} database \cite{PecanStreet}. In particular, we show that the equilibrium induced by the hourly billing mechanism \cite{baharlouei2013achieving} is robust to the level given to preferences.

This paper is organized as follows: \Cref{sec: context and energy consum} introduces the notion of consumers temporal preferences and defines the energy consumption game model. In \Cref{sec: sys costs social cost}, we define the main metrics of our study: we recall the definition of the Price of Anarchy and define the Price of Efficiency. In \Cref{sec: properties theo}, we give theoretical results and properties on the formulated model. We present explicit results on the equilibria in a simplified framework. Last, \Cref{sec: numerical exp} is devoted to numerical experiments on a realistic framework, where we simulate the equilibria among thirty Texan residential consumers in January, 2016.

\section{Context and Energy Consumption Game}
\label{sec: context and energy consum}
The model of this work falls within the class of DSM studies where the interaction of individual consumers is coordinated introducing an energy consumption game, as in \cite{mohsenian2010autonomous}. While all the proposed results could be applied to numerous operational frameworks, the one described here consists in the interaction between a provider and its set $\mathcal{N} = \{1,\cdots,N\}$ of consumers in a given day. As opposed to \cite{mohsenian2010autonomous} which does not show any preference for consumers and suppose that they are indifferent to any consumption schedule as soon as it satisfies their constraints, here we will focus on the integration of individual preferences of consumers into their objectives.

Indeed, consumers tend to have a ``natural" or preferred consumption profile, and asking them to deviate from it might be inconvenient or decrease their comfort. Individual utility functions have been previously used through different models. A common approach (see for instance \cite{fahrioglu2001using} and \cite{samadi2012advanced}) is to consider that a consumer's utility can be modeled as an increasing function of the total energy he receives. Here, on the contrary, we keep the assumption made in \cite{mohsenian2010autonomous} that consumers have flexible appliances that need a fixed quantity of energy per day, and this demand must be satisfied each day. However, we assume that consumers are not indifferent to the time they can use electricity and therefore use their appliances. 

\subsection{Introducing users temporal preferences}
From the provider's point of view, only the load profile asked by each user $n$, $(\ell^h_n)_{h\in \mathcal{H}}$ matters, where $\mathcal{H}$ is the discrete set of time periods considered. However, that may not be the case for users: for instance, one would like to charge an Electric Vehicle (EV) battery as soon as possible in case of unscheduled need (Plug-and-Charge), or one would like to turn on the heating system in a household at precise time periods (\cite{nguyen2014}), etc.
We denote user $n$'s \textit{preferred} or \textit{desirable} consumption profile by the vector $\bm{\hat{\ell}}_{na}=(\hat{\ell}_{na}^h)_h$ for his flexible appliance $a$. As a result, user $n$ would like to receive the power profile $(\hat{\ell}_n^h)_h:= (\sum_{a} \hat{\ell}_{na}^h)_h$ and, if he has no incentives to do otherwise, this profile will be his actual one. 
 Deviating from the profile $\hat{\bm{\ell}}_n$ decreases the comfort or utility of consumer $n$. To model this fact, we introduce the individual utility $u_n(\bm{\x}_n)$ of consumer $n$  as the opposite of the squared distance\footnote{
 In general, one could use $ u_n(\bm{\x}_n) := -\omega_{n}  d(\bm{\x}_n,\hat{\bm{\x}}_n) $ 
where $d(.,.)$ is an arbitrary metric. For simplicity and computational purposes, we use $d=\norm{.}^2$.} between the actual consumption profile of consumer $n$, $\bm{\x}_n$, and his preferred profile $\hat{\bm{\x}}_n$:
\vspace{-0.1cm}
\begin{equation}
\label{eq: quad utilities}
 u_n(\bm{\ell}_n) := -\omega_{n} \tssum_h ( \ell_n^h - {\hat{\ell}}_n^h)^2 \ ,
 \vspace{-0.1cm}
\end{equation}
where the weight $\omega_n$ indicates how much user $n$ values the distance to his preference. As some users will give more importance to their electricity bills (defined below) and some others to their utility, in a general framework we assume that each user could use a different weight $\omega_n$ in this model.

\subsection{Users billing mechanism}
\label{subsec: users billing}
As done in \cite{mohsenian2010autonomous,Paulin2017ISGT},  we suppose that the providing costs on each period $h$ are represented as a quadratic function $\tilde{C}(L^h)$ of the total load $L^h=\ell^h_{\text{NF}}+\ell^h$, where $\ell^h_{\text{NF}}$ denotes the aggregated nonflexible load at period $h$ and $\ell^h$ the flexible part:
\vspace{-0.1cm}
\begin{equation}
\label{eq: uniform quad system costs}
\tilde{C}(L^h):=\tilde{a}_0+ \tilde{a}_1 L^h + \tilde{a}_2 {(L^h)}
^2.
\vspace{-0.1cm}
\end{equation}
The surplus cost induced by the flexible part of the load $\ell^h$, at time $h$, denoted by $C_h(\ell^h)$, can be deduced as:
\begin{align}
\label{eq: cost flexible load}
C_h(\ell^h):&=\tilde{C}(L^h)-\tilde{C}(\ell^h_{\text{NF}}) = a_{1,h} \ell^h + a_2 ({\ell^h})^2  
\end{align}
with $a_{1,h}:=\left( \tilde{a}_1+2\tilde{a}_2 \ell^h_{\text{NF}} \right) $ and $a_2=\tilde{a}_2 $. Even if the cost function for the provider \eqref{eq: uniform quad system costs} does not depend on time, the nonflexible load profile $\ell^h_{\text{NF}}$ induces a difference in the  costs  $(C_h)_{h\in\mathcal{H}}$ between the different time periods $h\in \mathcal{H}$.

In this paper, we consider that the nonflexible part of the load is managed and billed  in a distinct process, e.g., in a standard contract. We focus on DR billing mechanisms for the \emph{flexible} part of the load. Through this study, we will consider two different billing mechanisms which, in practice, would require a ``two-way'' communication system \cite{ipakchi2009}, which enables the system operator to send its price functions $(C_h)_h$ and aggregated load $(\x^h)_h$ and users to send back their consumption profile $\bm{\x}_n$. 
  First, we consider the Daily Proportional (DP) billing mechanism introduced in \cite{mohsenian2010autonomous}: we assume that the system costs $C_h(\ell^h)$  induced by the flexible load at time $h$:
\vspace{-0.1cm}
\begin{equation}
\ell^h:= \textstyle\sum_n \ell^h_n \
\vspace{-0.1cm} 
\end{equation}
are shared among users proportionaly to their total flexible consumption on the entire day $E_n=\sum_{h\in\mathcal{H}} \ell_n^h$. Formally, each consumer will pay the daily bill:
\begin{equation}
\label{eq: bills DLP}
b_n^{\text{DP}}(\bm{\ell_n},\bm{\ell_{-n}}) = \frac{E_n}{E} \sum_{h\in\mathcal{H}} C_h(\ell^h) \ ,
\end{equation}
where $\bm{\ell_{-n}}=(\bm{\ell}_m)_{m\neq n}$ and $E:=\sum_n E_n$.
We will compare its efficiency to the natural ``congestion" Hourly Proportional (HP) billing mechanism introduced in \cite{baharlouei2013achieving}, where system costs on each period are shared among consumers respectively to their consumption on this period. Formally, the daily bill $b_n^{\text{HP}}$ of user $n$ for his flexible consumption is\footnote{Introducing  per-unit prices $c_h:= C_h(\ell^h)/ \ell^h$, the bill of $n$ can also be formulated in the ``congestion" form: $b_n^{\text{HP}}(\bm{\ell}) =\ell^h_n c_h({\ell}^h)$ analyzed in \cite{orda1993competitive}.}:
\begin{equation}
\label{eq: bills HLP}
b_n^\text{HP}(\bm{\ell_n},\bm{\ell_{-n}})= \sum_{h\in \mathcal{H} } \frac{\ell^h_n}{\ell^h} C_h(\ell^h) \ .
\end{equation}
Intuitively, with $b_n^\text{HP}$, users are more impacted by their actions to use expensive/cheap time periods than with $b_n^\text{DP}$ where the costs induced by actions are ``averaged'' over the day. This property helps to interpret the results of \Cref{sec: properties theo} and \Cref{sec: numerical exp}.

\subsection{Energy consumption game}
To analyze the impact of the importance given to users' temporal preferences, we consider through this work the parametrized users' objective functions:
\begin{equation}
\label{eq: obj n param}
f_n^\alpha(\bm{\ell_n},\bm{\ell_{-n}}):= (1-\alpha) b_n( \bm{\ell} ) - \alpha u_n(\bm{\ell})  
\end{equation}
where the \emph{preference factor} $\alpha\in [0,1]$ indicates the weight given to user $n$'s preference\footnote{We could extend this study by using different $(\alpha_n)_n $ for different users.} in comparison to his bill $b_n$.  We get the following optimization problem for user $n$:
\begin{subequations}
\label{eq:user_problem}
\begin{align}
&\min_{\bm{\ell}_n \in \mathbb{R}^H }  f_n^{\alpha}(\bm{\ell}_n,\bm{\ell}_{-n}) \\
\label{cons: HP total power} \text{s.t. } & \hspace{0.5cm}   \textstyle\sum_{h\in \mathcal{H}} \ell^h_{n} = E_{n}  , \\ 
 \label{cons: HP minmax power}& \hspace{0.5cm} \underline{\x}^h_{n} \leq \x_{n}^h 
\leq \overline{\x}^h_{n} , \forall h \in \mathcal{H} \ .
  \end{align}
\end{subequations}
Constraint (\ref{cons: HP total power}) expresses that a fixed daily amount of energy is required for the flexible appliances of user $n$ (EV battery, washing machine...). Due to physical limits of his electrical items or personal constraints, the power given to $n$ is bounded (\ref{cons: HP minmax power}). 
We denote more compactly by $\mathcal{L}_n$ the feasible set of user $n$, given by the polytope \eqref{cons: HP total power}, \eqref{cons: HP minmax power}, and $\mathcal{L}:=\mathcal{L}_1 \times \dots \times \mathcal{L}_N$.
When $\alpha=0$, user's preference $\hat{\bm{\ell}}_n$ has no influence on his behaviour, while when $\alpha=1$, the user gives no importance to $b_n$ and only wants to minimize $-u_n$: his resulting load profile will be exactly his preference $\hat{\bm{\ell}}_n$.

As $f_n^\alpha$ depends on the consumption of $n$ but also on other users, this induces a game between users \cite{fudenberg1991game} denoted by $\mathcal{G}_\alpha := \left( \mathcal{N},\mathcal{L},(f_n^\alpha)_{n}  \right)$. We will use the notations $\mathcal{G}_\alpha^\text{DP}$ and $\mathcal{G}_\alpha^\text{HP}$ when we specify the billing mechanism according to the DP rule \eqref{eq: bills DLP} or HP rule \eqref{eq: bills HLP}. 
The importance that each user gives to his utility function $u_n$ in comparison to his bill $b_n$, through the parameter $\alpha$, will change the set of Nash Equilibria (NE) of the game $\mathcal{G}_\alpha$ given by:
\begin{equation*}
\mathcal{L}^{\text{NE}}_\alpha := \{ \bm{\ell}\in \mathcal{L}: \forall n, \forall \bm{\ell}'_n \in\mathcal{L}_n, \  f_n^{\alpha}(\bm{\ell}_n,\bm{\ell}_{-n}) \leq f_n^{\alpha}(\bm{\ell}'_n,\bm{\ell}_{-n})\}.
\end{equation*}
 
\section{Social cost versus System costs}
\label{sec: sys costs social cost}

An Independent System Operator is interested in both an efficient electricity network and the welfare of the consumers. Starting with the latter, we define the \emph{social cost} of the set of consumers as the sum of their objective functions:
\begin{equation}
\label{eq:social cost}
\hspace{-0.0cm}\SC_\alpha(\bm{\ell})\hspace{-0.1cm}:=\hspace{-0.1cm}\sum_{n\in\mathcal{N}} \hspace{-0.1cm} f_n^\alpha(\bm{\ell}) \hspace{-0.1cm} =  (1-\alpha) \hspace{-0.1cm} \sum_{n\in\mathcal{N}} b_n(\bm{\ell}) - \alpha \hspace{-0.1cm} \sum_{n\in\mathcal{N}}u_n(\bm{\ell}_n) \ .
\end{equation}
To quantify the efficiency of a billing mechanism in a game, we consider the standard notion of Price of Anarchy (PoA) introduced by Koutsoupias and Papadimitriou. The PoA measures the gap between the minimal social cost \eqref{eq:social cost}, and the social cost induced by the worst equilibrium of the game.
\begin{definition}[\cite{koutsoupias1999worst}]{Price of Anarchy (PoA). \\}
\label{def:PoA}
Given a game $ \mathcal{G}$, $\mathcal{X}_{\mathcal{G}}^{\text{NE}}$ its set of Nash Equilibria 
and $\SC^*
$ its minimal social cost, the price of anarchy of $\mathcal{G}$ is given as:
\begin{equation}
\label{eq:PoA}
\text{PoA}(\mathcal{G}):=\left(	 {\textstyle\sup_{\bm{\ell}\in \mathcal{X}_{\mathcal{G}}^{\text{NE}}} \SC\left(\bm{\ell}\right)} \right)  /  \  {\SC^*} \ .
\end{equation}
\end{definition}

From the provider's point of view, only the costs induced for the system, without users personal utilities, matters: we denote by $\mathcal{C}$ the total system costs function, defined as the providing costs induced by a consumption profile $(\ell^h)_{h\in\mathcal{H}}$:
\begin{equation}
\label{eq:system cost}
\mathcal{C} (\bm{\ell}):= \textstyle\sum_{h\in \mathcal{H}} C_h(\ell^h) \ .
\end{equation}
Note that in the particular billing mechanisms considered in \eqref{eq: bills DLP} and \eqref{eq: bills HLP}, we assume that the provider costs are shared among users\footnote{This assumption could be relaxed, as done in \cite{mohsenian2010autonomous}, by adding a ratio profit $\kappa>1$ for the provider, so that we have the equality $\sum_n b_n(\bm{\ell})= \kappa \mathcal{C} (\bm{\ell})$.}, so that we have the equality $\mathcal{C} (\bm{\ell}) = \sum_n b_n(\bm{\ell}) $.

We introduce a measure similar to the PoA \eqref{eq:PoA}, but that will be more relevant for a provider that is more interested in the system costs $\mathcal{C}$ and does not have access to the utility functions $(u_n)_n$ of its users:

\begin{definition}{Price of Efficiency (PoE).}
\label{def:PoE}

Given a game $ \mathcal{G}$, $\mathcal{X}_{\mathcal{G}}^{\text{NE}}$ its set of Nash Equilibria 
 and its minimal feasible system costs $\mathcal{C}^*:= \min_{\bm{\x} \in \mathcal{L}} \mathcal{C}(\bm{\x})$ , the price of efficiency of $\mathcal{G}$ is given as:
\begin{equation}
\label{eq:PoE}
\text{PoE}(\mathcal{G}):=\left( { \textstyle\sup_{\bm{\ell}\in \mathcal{X}_{\mathcal{G}}^{\text{NE}}} \mathcal{C}\left(\bm{\ell}\right)} \right) / \ {\mathcal{C}^*} \ .
\end{equation}
\end{definition}
Observe that $\text{PoA} \geq 1$ and $\text{PoE} \geq 1$. Following \eqref{eq: obj n param}, one can notice that for $\alpha=0$, $\text{PoE}(\mathcal{G}_\alpha)=\text{PoA}(\mathcal{G}_\alpha)$. In general, the PoA and PoE will be different as shown below.
%
%

\section{Properties}
\label{sec: properties theo}
\subsection{Potential property, existence and uniqueness of NE}
We start by showing that the considered games have the property of potential (see \cite{monderer1996potential}).
\begin{theorem}
\label{th: DLP potential}
$\mathcal{G}^\text{DP}_\alpha$ is a weighted potential game with potential:
\begin{equation}
\vspace{-0.1cm}
W^\text{DP}_\alpha= (1-\alpha) \sum_{h\in\mathcal{H}} C_h(\ell^h) - \alpha \sum_n \frac{E}{E_n} u_n(\bm{\ell}_n) \ .
\vspace{-0.1cm}
\end{equation}
\textit{Proof:} $ \forall n$, $\nabla_nf_n^\alpha=\frac{E_n}{E}\nabla_n W^\text{DP}_\alpha$, we conclude from \cite{monderer1996potential}.
\end{theorem}
With the billing mechanism HP, we need an additional assumption on the system cost functions to get a similar result: 
\begin{theorem}
\label{th:potential}
If we consider quadratic costs \eqref{eq: uniform quad system costs}, $\mathcal{G}_{\alpha}^\text{HP}$ is an exact potential game with potential:
\begin{equation}
\label{eq: potential def}
\hspace{-0.1cm}W^\text{HP}_\alpha\hspace{-0.1cm}= \hspace{-0.1cm} (1-\alpha) \hspace{-0.1cm} \left[ \hspace{-0.05cm}\sum_{h\in \mathcal{H}} \hspace{-0.1cm}\textstyle\frac{a_{2,h}}{2} \hspace{-0.05cm}\Big( \hspace{-0.1cm}{(\ell^h)}^2 \hspace{-0.1cm} +\hspace{-0.1cm} \sum_{n} \hspace{-0.1cm} {(\ell^h_n)}^2 \Big)\hspace{-0.1cm} + \hspace{-0.1cm} a_{1,h} \ell^h \right] \hspace{-0.1cm} - \alpha \hspace{-0.1cm}\sum_{n\in\mathcal{N}} \hspace{-0.1cm} u_n(\bm{\ell}) \ .
\end{equation}
\textit{Proof:} Similarly to the proof of \Cref{th: DLP potential}, $ \forall n$,  $\nabla_nf_n^\alpha=\nabla_nW^\text{HP}_\alpha$.
\end{theorem}

From the fact that $W^\text{DP}_\alpha$ and $W^\text{HP}_\alpha$ are strictly convex and from \cite{monderer1996potential}, we can deduce the existence and uniqueness of NE:
\begin{Cor}
\label{cor: unicity exist NE}
In the games $\mathcal{G}^\text{DP}_\alpha$  and $\mathcal{G}^\text{HP}_\alpha$ there exists a unique Nash Equilibrium corresponding respectively to the minimum argument of $W^\text{DP}_\alpha$  and of $W^\text{HP}_\alpha$ over the set $\mathcal{L}$.
 \end{Cor} 
\Cref{cor: unicity exist NE} extends the results of \cite{mohsenian2010autonomous} which gives the uniqueness of NE in the particular case of $\alpha\hspace{-0.1cm}=\hspace{-0.1cm}0$ with the DP billing.
 
A natural algorithm to compute a NE is to run the Best Response Dynamics (BRD), as defined below.
\begin{definition}[\cite{gilboa1991social}]{Best Response Dynamics (BRD).\\}
\label{def:BRD}
At each iteration $k$, a user $n_k$ is randomly chosen and solves problem \eqref{eq:user_problem} to optimum $\bm{\ell}^{*}_{n_k}$, with load of others $\x_{-n_k}^{(k)}$ fixed ($n_k$ best responses to the others). We update $\bm{\ell}^{(k+1)}_{n_k}=\bm{\ell}^{*}_{n_k}$.
\end{definition}
From \Cref{th: DLP potential} and \Cref{th:potential}, we deduce the convergence of BRD:
\begin{Cor}
In $\mathcal{G}^\text{DP}_\alpha$  and $\mathcal{G}^\text{HP}_\alpha$, BRD is equivalent to a block coordinate minimization of the potential function. Hence, it converges to the unique NE of the game (see \cite{beck2013convergence}).
\end{Cor}
%
%

%
%
\subsection{Theoretical results on a simplified framework}
\label{subsec: theo results toy}

In this section, we consider that the set $\mathcal{H}$ is reduced to two time periods $\mathcal{H}:=\{P,O\}$ which represent for instance the \emph{Peak} and \emph{Offpeak} times. For computational purposes, we consider that the system costs are reduced to a quadratic term:
\begin{equation}
\label{eq: toy quad costs}
\forall h \in \mathcal{H}, \ C_h= ({\ell^h})^2 \ ,
\end{equation}
and there is no nonflexible part as in the general case described in \Cref{subsec: users billing}.  
Each consumer $n$ has a preference weight $\omega_n=1$, a preferred profile $(\hat{\ell}_n^P,\hat{\ell}_n^O)$, satisfying $\hat{\ell}_n^P+\hat{\ell}_n^O=E_n$ as in \eqref{cons: HP total power}. Without loss of generality, we assume $\hat{\ell}^P \geq \frac{E}{2} \geq \hat{\ell}^O $.
Power constraints \eqref{cons: HP minmax power} are replaced by positivity $\x^h_n \geq 0$.

\subsubsection{Nash equilibrium}

From the KKT conditions of optimality, we get the following result:
\begin{theorem}
Assume that for all $n\in \mathcal{N}$, we have:
\begin{equation}
\label{eq:pos DLP assum}
\frac{\hat{\x}^P_n}{E_n}+\frac{1}{2}\geq \frac{\hat{\x}^P}{E} \ ,
\end{equation}
then, for $ \alpha \in (0,1]$, the unique NE of $\mathcal{G}_\alpha^\text{DP}$ is given by: 
\begin{align}
\label{eq:l^h expr DLP toy}
\ell^P_n&= \hat{\ell}^P_n + \tsfrac{E_n}{E} \tsfrac{1-\alpha}{2}(\hat{\x}^{O}-\hat{\x}^P) \ ,
\end{align}
with symmetric expression for $\x^O_n$. For $\alpha\hspace{-0.1cm}=\hspace{-0.1cm}0$, the KKT system is degenerated, and any $(\bm{\x}^h_n)_n$ satisfying \eqref{eq:l^h agg expr DLP toy} below is a NE.

Assume that for all $n \in \mathcal{N}$, we have:
\begin{equation}
\label{eq:pos HLP assum}
{2(N-1)} \hat{\x}^P_n  \geq {(\hat{\x}^P-\hat{\x}^O)-E_n}\ ,
\end{equation}
then, for $ \alpha \in [0,1]$, the unique NE of $\mathcal{G}_\alpha^\text{HP}$ is given by: 
\begin{align}
\label{eq:l^h expr HLP toy}
\ell^P_n&=  \hat{\ell}^P_n + \textstyle\frac{1-\alpha}{2(1+\alpha)}\left( \phi( \alpha)(\hat{\x}^{O}-\hat{\x}^P) + (\hat{\ell}^{O}_n-\hat{\ell}^P_n) \right), 
\end{align}
with symmetric expression holding  for $\x^O_n$, and with:
\begin{equation}
\phi(\alpha)\eqDef \textstyle\frac{2\alpha}{(1+\alpha)+(1-\alpha)N} \in[0,1] \ .
\end{equation}
\end{theorem}

 One can check that the positivity of the offpeak load $\x^O_n$ in \eqref{eq:l^h expr DLP toy} and in \eqref{eq:l^h expr HLP toy} is always verified. The positivity of the peak load $\x^P_n$ is a consequence of assumptions \eqref{eq:pos DLP assum} and \eqref{eq:pos HLP assum}.
 
 We consider that  \eqref{eq:pos DLP assum} and \eqref{eq:pos HLP assum} hold through all this Section.


\begin{Cor}
The aggregated load at the NE is given by:
\begin{align}
\label{eq:l^h agg expr DLP toy}
\text{ for } \mathcal{G}_\alpha^\text{DP} ,  \ \ell^P&= \textstyle\frac{E}{2}+\alpha \textstyle\frac{(\hat{\ell}^P-\hat{\ell}^{O} ) }{2} \ , \\
\label{eq:l^h agg expr HLP toy}
\text{ for } \mathcal{G}_\alpha^\text{HP} ,  \ \ell^P &= \textstyle\frac{E}{2}+ \phi(\alpha)\textstyle\frac{(\hat{\ell}^P- \hat{\ell}^{O})}{2} \ .
\end{align}

\end{Cor}

 With HP and DP, the aggregated load evolves to the preferred profile when $\alpha$ goes to one but, with the HP mechanism, this evolution is influenced by the number of players $N$.
 
\subsubsection{System Costs}

The total system costs $\mathcal{C}=(\ell^P)^2+ (\x^O)^2$ at the equilibrium are given from \eqref{eq:l^h agg expr DLP toy} and \eqref{eq:l^h agg expr HLP toy} by:
\begin{align}
\label{eq : System costs alpha DP}
\text{ for } \mathcal{G}_\alpha^\text{DP} ,  \ &\mathcal{C}^\text{DP} (\alpha):= \textstyle\frac{1}{2} \left( {E^2}+ {\alpha^2 (\hat{\ell}^P- \hat{\ell}^O)^2}\right) \ , \\
\label{eq : System costs alpha HP}
\text{ for } \mathcal{G}_\alpha^\text{HP} ,  \ &\mathcal{C}^\text{HP}(\alpha):= \textstyle\frac{1}{2} \left({E^2}+ {\phi^2(\alpha)(\hat{\ell}^P- \hat{\ell}^O)^2} \right) \ .
\end{align}

\begin{figure}[!t]
\vspace{-0.5cm}
\centering
\includegraphics[width=0.42\textwidth]{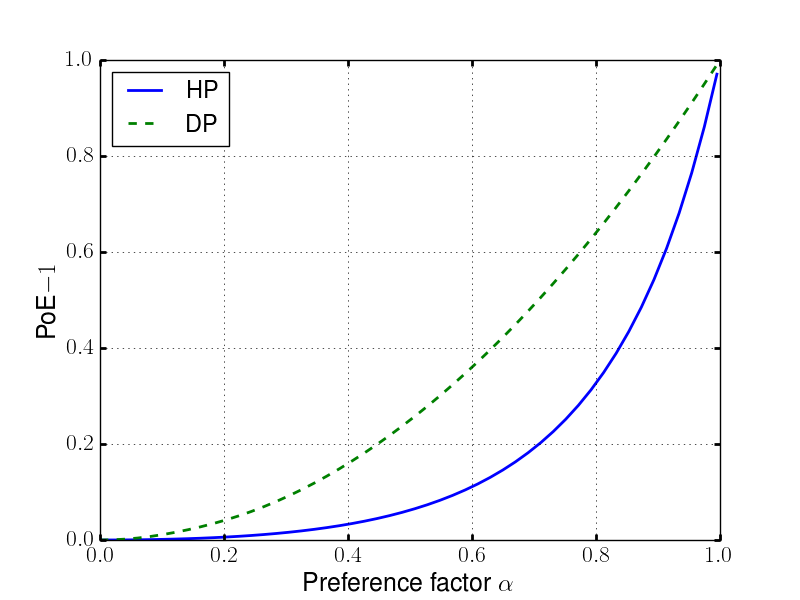}
\caption{Evolution of PoE-1 with costs \eqref{eq: toy quad costs} and $N=5$ users. \\ \textit{For all $\alpha$, HP billing is more efficient for the system than DP.}}
\label{fig: evo PoA PoE alpha}
\vspace{-0.3cm}
\end{figure}

On \Cref{fig: evo PoA PoE alpha}, we see that the $\text{PoE}$  is increasing with $\alpha$ in both cases (the proof is straightforward from \eqref{eq : System costs alpha DP} and \eqref{eq : System costs alpha HP}), and that it is always smaller with the HP billing, as shown below:

\begin{theorem}
The system costs induced by the equilibrium with HP are always smaller than with DP, or equivalently:
\begin{equation}
\forall \alpha \in [0,1], \ \text{PoE}(\mathcal{G}_\alpha^\text{HP} ) \leq \text{PoE}(\mathcal{G}_\alpha^\text{DP} ) 
\end{equation}
and the inequality is strict for $\alpha\in (0,1)$.
\end{theorem}
\textit{Proof:} First, note that in \eqref{eq:PoE}, the minimal system cost $\mathcal{C}^*$ does not depend on $\alpha$, so that $\mathcal{C}$ and PoE are proportional. From the expressions \eqref{eq : System costs alpha DP} and \eqref{eq : System costs alpha HP}, we get that: \\
$\mathcal{C}^\text{DP}(\alpha)-\mathcal{C}^\text{HP}(\alpha) =  \frac{\alpha^2 (\hat{\ell}^P- \hat{\ell}^O)^2}{2} \left(1- \textstyle\frac{4}{(N(1-\alpha) + (1+\alpha))^2} \right)
$
$ >0$ because $ \forall \alpha \in (0,1), \frac{4}{(N(1-\alpha) + (1+\alpha))^2}<1 \Longleftrightarrow N>1$.

\Cref{fig: evo PoA PoE alpha} shows the  evolution of the PoE induced by the NE of $\mathcal{G}_\alpha$, in the case of $N=5$ players that have a flexible energy $E_n=1$ that they prefer loading totally on peak hour ($\hat{\ell}^P_n=1)$.

\subsubsection{Social Cost}

\begin{figure}[!t]
\centering
\hspace{-0.5cm}
\subfloat{\includegraphics[width=0.225\textwidth]{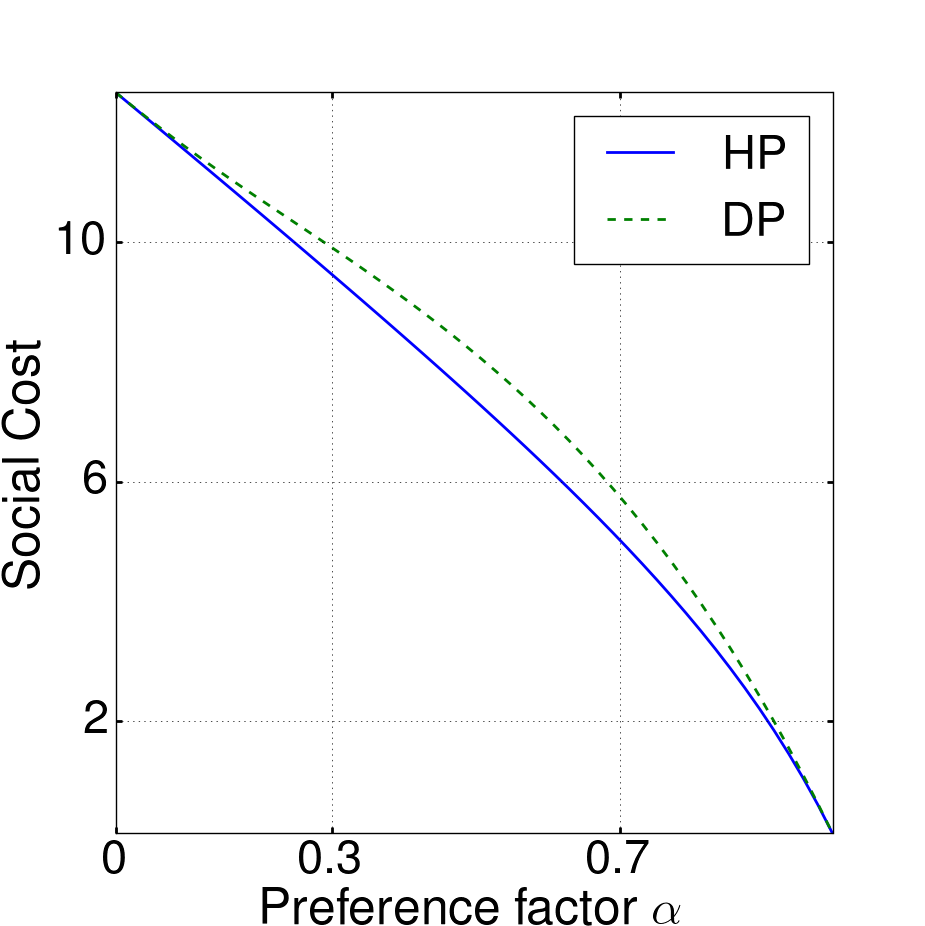}
\label{fig: evo toy social cost alpha}
} \subfloat{
\includegraphics[width=0.29\textwidth]{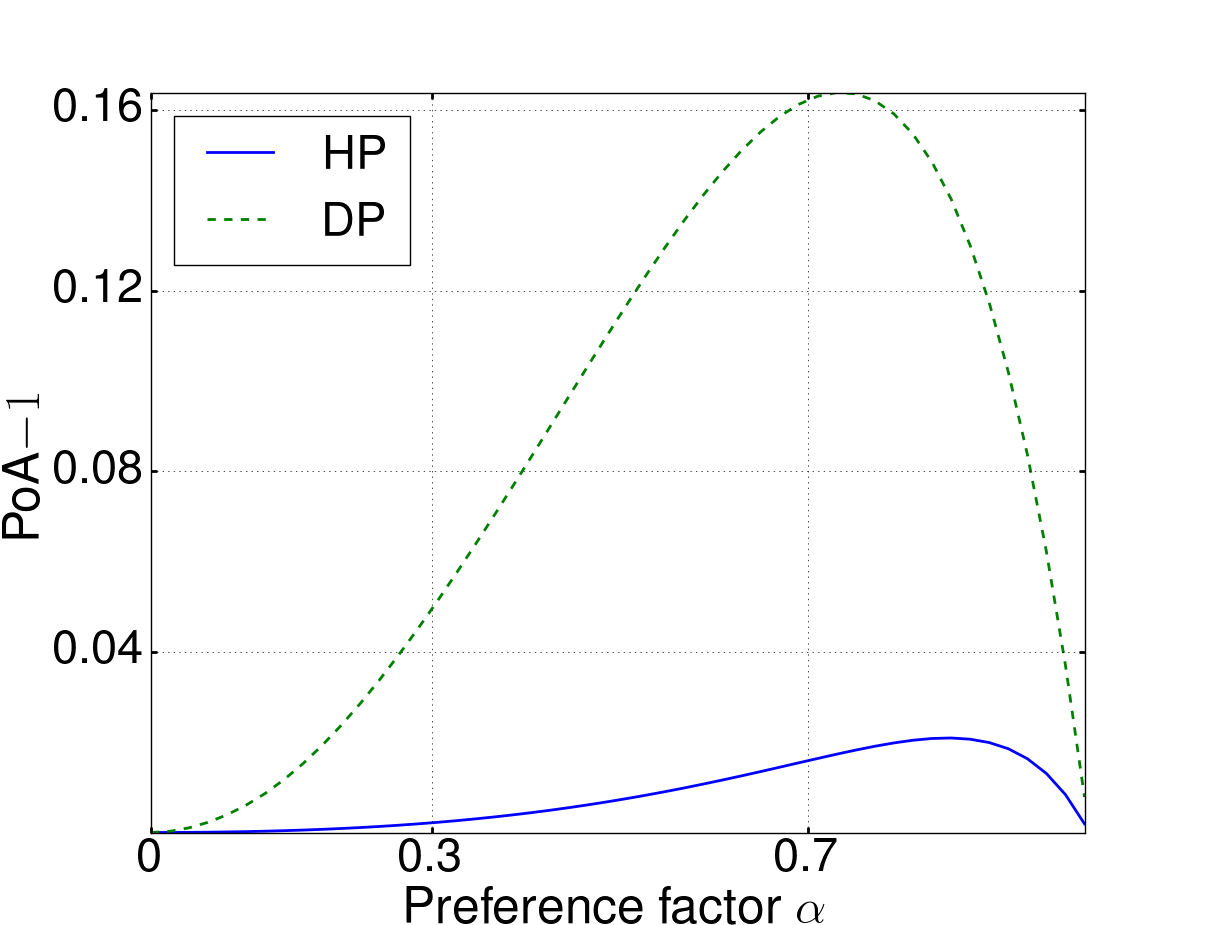}
}
\caption{Evolution of Social Cost and PoA-1 with costs  \eqref{eq:  toy quad costs}.\\ \textit{Social Cost always decreases with $\alpha$. The PoA is unimodal and with DP, it reaches its maximum at a critical level of $\alpha \simeq 0.72$, where it is more than 10\% larger than with HP. }}
\vspace{-0.3cm}
\end{figure}

If users do not care about their bills but only on their utility ($\alpha=1$), they choose their preferred profile $(\hat{\ell}^P_n,\hat{\ell}^O_n)$. As a result, the social cost will be exactly zero. On the opposite, if consumers only care about their bills ($\alpha=0$), \cite{mohsenian2010autonomous} shows that users will reach the optimal system cost in $\mathcal{G}_0^\text{DP}$ (the potential $W_0^\text{DP}$ is equal to the system costs) while \cite{Paulin2017ISGT} shows that the equilibrium in $\mathcal{G}_0^\text{HP}$ will stay close to the social optimum (it is even optimal in the framework of this section, as seen in \eqref{eq:l^h agg expr HLP toy}). However, it is not clear how the social cost evolves with $\alpha \in [0,1]$. \Cref{fig: evo toy social cost alpha} shows that with both the DP and HP mechanism, the social cost is a decreasing function of $\alpha$. We prove this in \Cref{th:SC DP decreasing} for the DP mechanism.

To this end, considering the expressions of the equilibrium in $\mathcal{G}_\alpha^\text{DP}$  from \eqref{eq:l^h expr DLP toy}, we get the induced social cost:
\begin{align}
\label{eq: SC DLP equi}
\SC^{\text{DP}}_\alpha&=(1-\alpha)\left[\textstyle\frac{E_{}^2}{2}+ \textstyle\frac{D^2}{2}(\alpha^2 + V_E (1-\alpha) \alpha ) \right]
\end{align}
with $E \eqDef 
 \sum_{n\in\mathcal{N}}E_n$, $D \eqDef (\hat{\x}^P-\hat{\x}^O )$ and $V_E \eqDef \sum_n \frac{E_n^2}{E^2}$.

\begin{theorem}
\label{th:SC DP decreasing}
$\SC^{\text{DP}}_\alpha$ is a decreasing function of $\alpha$.
\end{theorem}
\textit{Proof.} $\partial_{\alpha} {\SC^\text{DP}_\alpha}\hspace{-0.1cm}=\hspace{-0.1cm}\frac{D^2}{2} \left[ -3(1\hspace{-0.1cm}-\hspace{-0.1cm}V_E)\alpha^2\hspace{-0.1cm}+\hspace{-0.1cm}2(1\hspace{-0.1cm}-\hspace{-0.1cm}2V_E)\alpha   \right] \hspace{-0.1cm}+\hspace{-0.1cm} D^2V_E\hspace{-0.1cm}-\hspace{-0.1cm}E^2 $ is always negative. Details are omitted here for brevity.
We did not manage to prove a symmetric result for 
$\SC^{\text{DP}}_\alpha$.

\section{Numerical experiments}
\label{sec: numerical exp}

In this Section, we present numerical results on the sensibility of the equilibria of $\mathcal{G}_\alpha^\text{DP}$ and $\mathcal{G}_\alpha^\text{HP}$ to $\alpha$, in a realistic framework. We simulate the games $\mathcal{G}_\alpha^\text{DP}$ and $\mathcal{G}_\alpha^\text{HP}$ and the convergence to the equilibria day by day on the set of the thirty one days of January 2016,  which we denote by $\mathcal{D}$, each day being decomposed by a hourly timeset $\mathcal{H}=\{0,1,\dots,23\}$.

\subsection{Parameters}
\paragraph{Consumers} We extracted $N=30$ residential consumption profiles of Electric Vehicles (EV) owners from \textit{PecanStreet Inc.} \cite{PecanStreet}, a database of residential consumers in Texas (U.S.). Each consumer has a nonflexible consumption $(\ell_{\text{NF},n}^h)_h$ (lights, cooking, TV...) and we consider EV charging as the flexible usage. We take the EV historical profile of user $n$ as its preferred profile $ (\hat{\x}^h_n)_h$, and assume it corresponds to its flexible energy need $E_n:=\sum_h \hat{\x}^h_n$. For power constraints \eqref{cons: HP minmax power}, we take $\underline{\x}^h_n=0$ and $\overline{\x}_n^h$ equal to the max observed value if hour $h$ was ever used by $n$ and $\overline{\x}_n^h=0$ otherwise.

\paragraph{System Costs}
As explained in \Cref{subsec: users billing}, we suppose that system costs are, for each time $h$, function of the total load $L_h=\ell_{\text{NF}}^h + \ell^h$, and are given in dollar cents as $\tilde{C}(L^h)=71.1 -4.17 L^h + 0.295 (L^h)^2 $. To compute those coefficients, we make an interpolation based on three load values and three corresponding prices ($C_h(L^h)/L^h$). The three  load values are the  mean ($33.8$kW), min ($17.8$kW) and max ($58.9$kW) values of the nonflexible load per hour aggregated  over the set $\mathcal{N} $ of consumers in all hours of January, 2016. 
The three corresponding prices are those proposed by the Texan provider Coserv (\cite{Coserv}): 8.0\textcent/kWh for base contracts, 14.0\textcent/kWh (Peak) and 5.5\textcent/kWh (OffPeak) in Time-of-Use contracts. From \eqref{eq: cost flexible load}, the cost of flexible load is given by:
\begin{equation}
C_h(\ell^h):= (-4.17+ 0.590\ell_\text{NF}^h) \ell^h + 0.295 (\ell^h)^2 \ . 
\end{equation}

To ensure that $b_n$ and $u_n$ are of the same order of magnitude, we use a common factor in \eqref{eq: quad utilities} of $\omega_n\hspace{-0.1cm}=\hspace{-0.1cm}\omega\hspace{-0.1cm}:=\hspace{-0.1cm}\frac{\mathcal{C}^*}{ {\sum_n \norm{ \bm{\ell}^*_n-\hat{\bm{\ell}}_n }^2} }\hspace{-0.1cm}=\hspace{-0.1cm}49.1$\textcent/kWh$^2$ with $\mathcal{C}^*\hspace{-0.05cm}=\hspace{-0.05cm}\mathcal{C}(\x^*)\hspace{-0.1cm}=\hspace{-0.1cm}\min_{\x} \mathcal{C}(\ell) $ the optimal costs.

Note that for $\alpha=1$, $\SC^*_\alpha= 0 =\SC^\text{DP}_{\alpha}=\SC^\text{HP}_{\alpha}$ so the PoA is not defined, but \Cref{fig: evo PoA real case} shows that $\lim_{\alpha \rightarrow 1} \text{PoA}( \mathcal{G}_\alpha) = 1$.

\subsection{Results}
For both mechanisms HP and DP, we compute the NE by playing a BRD (\Cref{def:BRD}) with a limiting number of 150 BR iterations, which in practice was sufficient for convergence.

The optimization problem \eqref{eq:user_problem} is a quadratic program that we solve with the optimization solver Cplex 12.6. Playing the BRD takes around 2.5sec. for each of the 50 values of $\alpha$ and each day in $\mathcal{D}$. The total simulation time was  3160 sec. with an Intel Xeon CPU E3-1240v3@3.4GHz$\times$8 run on 5 threads.

\begin{figure}[!t]
\centering
\includegraphics[width=0.495\textwidth]{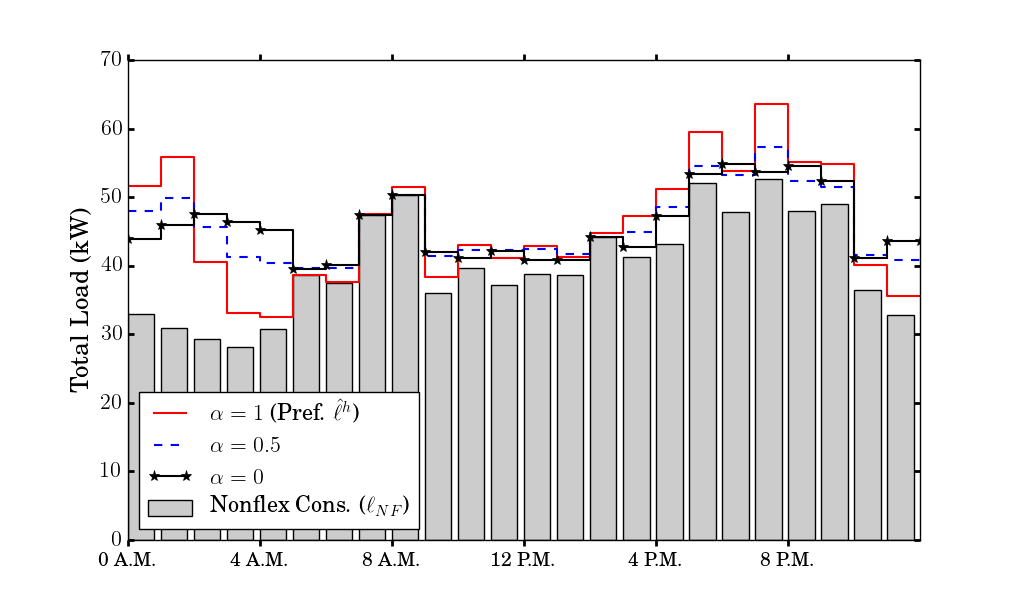}
\caption{Aggregated consumption profiles at equilibrium $\mathcal{G}^\text{HP}_\alpha$ for $\alpha=0, \ 0.5 $  and 1 on day 10/01/2016 with 30 users. \\
\textit{The equilibrium profile converges to $(\hat{ \ell}^h)_h$ when $\alpha \rightarrow 1$.}  }
\label{fig: evo profile alpha}
\vspace{-0.3cm}
\end{figure}

\Cref{fig: evo profile alpha} shows the different aggregated profile $(\ell^h)_h=(\sum_n \ell^h_n)_h$ at the equilibrium of $\mathcal{G}_\alpha^\text{HP}$ on January, $10$, chosen arbitrarily in $\mathcal{D}$. We can see a significant variation  (more than $15\%$) on the aggregated load when $\alpha$ changes.

\begin{figure}[!t]
\vspace{-0.4cm}
\centering
\includegraphics[width=0.48\textwidth]{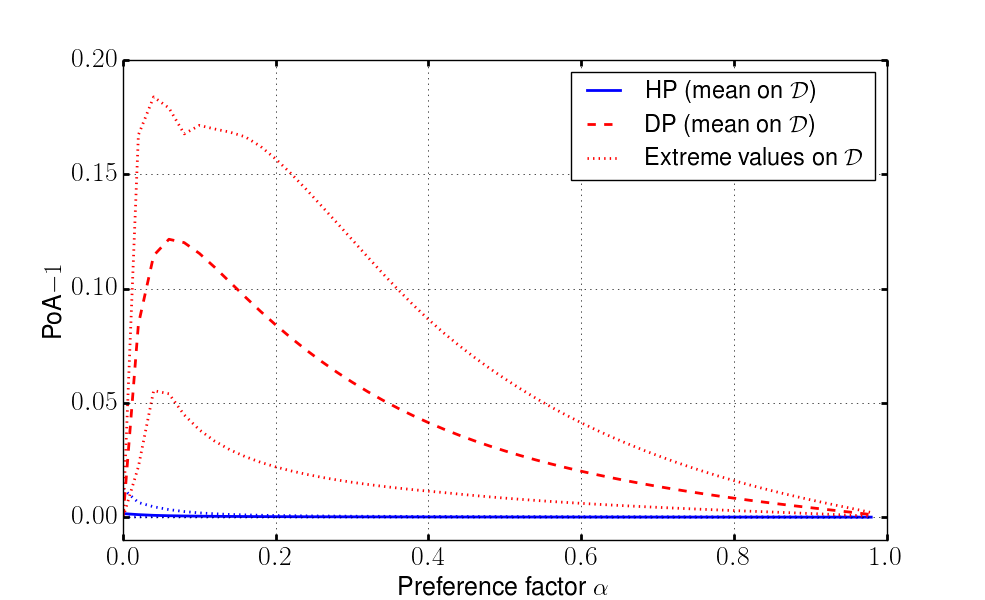}
\caption{Evolution of PoA-1 (mean on days) with $\alpha$. \\\textit{HP is more robust and has a smaller PoA than DP.}} 
\label{fig: evo PoA real case}
\end{figure}

\begin{figure}[!t]
\vspace{-0.54cm}
\centering
\includegraphics[width=0.48\textwidth]{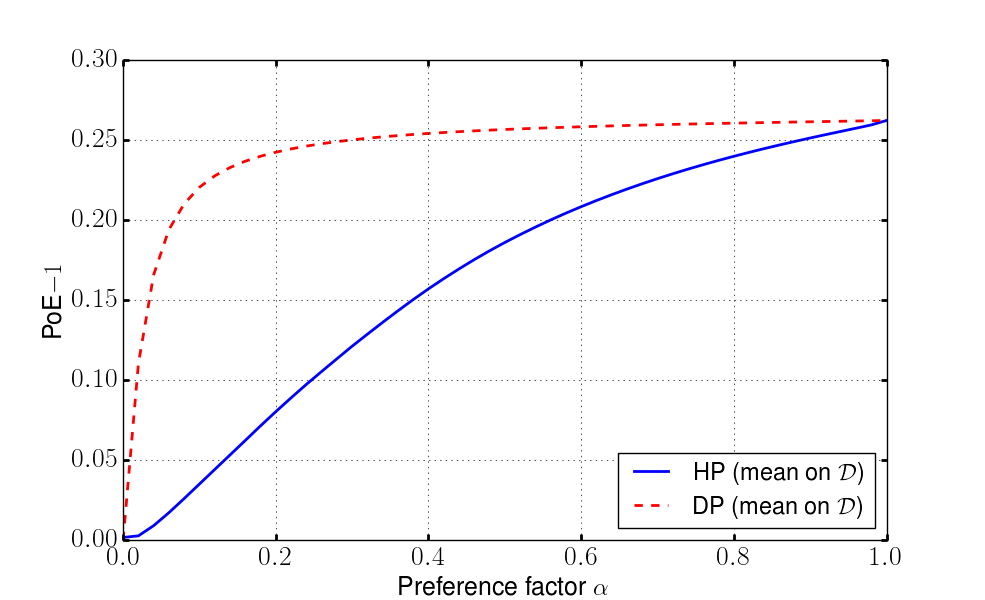}
\caption{Evolution of PoE-1 (mean on days) with $\alpha$. \\ \textit{For $\alpha >3\cdot 10^{-4}$, the equilibrium induced by the HP billing is also more efficient in terms of system costs}.} 
\label{fig: evo PoE real case}
\vspace{-0.3cm}
\end{figure}

From \Cref{fig: evo PoA real case}, we see that at $\alpha=0$, there is a very small PoA for the HP mechanism (see \cite{Paulin2017ISGT} for a deeper analysis) while DP achieves optimality. However, when $\alpha$ grows the HP mechanism becomes much more efficient than DP in terms of PoA (\Cref{fig: evo PoA real case}) and PoE (\Cref{fig: evo PoE real case}), as already seen in the simplified framework of \Cref{subsec: theo results toy}. We observe that, in \Cref{subsec: theo results toy} as in this realistic case,  PoA$(\mathcal{G}^\text{DP}_\alpha)$ is an unimodal function of $\alpha$.

\Cref{fig: evo PoA real case} shows that the PoA induced by the equilibrium of $\mathcal{G}_\alpha^\text{HP}$ remains very low (it is maximal at $\alpha=0$ with PoA$=$1.0015 and then decreases) while the PoA of DP reaches a maximum of 1.122 at $\alpha$=0.06: this billing mechanism is much less robust to consumers' preferences. This lack of robustness is underlined by the important discrepancy between the minimal and maximal PoA values over our set of 31 days.
\Cref{fig: evo PoE real case} shows that the PoE, as a function of $\alpha$, is much more concave for the DP mechanism, resulting in larger system costs on a wide range of $\alpha$ (the two curves intersect at $\alpha \simeq 3\cdot10^{-4}$). As a result, the HP mechanism will also be more interesting for the provider. For $\alpha \leq 3\cdot10^{-4}$ the system costs are bigger for HP than DP because of the small PoA mentioned before \cite{Paulin2017ISGT}.

\section{Conclusion}
We considered a game theoretic model to study the behavior of residential consumers in a DR program. We formulated an energy consumption game with a temporal preference term in each user's cost function. We gave several theoretical results on a simplified test case and showed by simulations that those results still hold in a realistic framework where consumers have a nonflexible load. Without consumers preferences, the Daily Proportional billing reaches the optimal social cost and is more efficient than the Hourly Proportional billing which is not exactly optimal. When we add the temporal preference term, the Hourly Proportional billing becomes much more advantageous than the Daily Proportional mechanism in terms of social cost and in terms of costs induced for the provider.

Several extensions of this work could be considered. First, the theoretical results could be extended to take into account a nonflexible part, or considering general functions instead of a quadratic model for system costs. Besides, we could study a dynamic population of users who have the choice to remain in the demand response program or not: if consumers are not satisfied with the program, they might consider another kind of contract or suscribing a more competitive provider.


%
%
%



\bibliographystyle{IEEEtran}
\bibliography{../../DRAFTS/Biblio_complete/biblio1,../../DRAFTS/Biblio_complete/biblio2,../../DRAFTS/Biblio_complete/biblioBooks}

\end{document}